\documentclass[11pt]{amsart}
\usepackage{amssymb,hyperref}
\usepackage{graphicx,float}
\usepackage{graphics} 
\usepackage{tikz}
\usetikzlibrary{arrows}

\newtheorem{theorem}{Theorem}

\newtheorem{conjecture}[theorem]{Conjecture}
\newtheorem{corollary}[theorem]{Corollary}
\newtheorem{definition}[theorem]{Definition}
\newtheorem{example}[theorem]{Example}

\newtheorem{exercise}[theorem]{Exercise}
\newtheorem{fact}[theorem]{Fact}
\newtheorem{lemma}[theorem]{Lemma}
\newtheorem{problem}[theorem]{Problem}
\newtheorem{proposition}[theorem]{Proposition}
\newtheorem{question}[theorem]{Question}
\newtheorem{remark}[theorem]{Remark}

\newcommand{\bcon}{\begin{conjecture}}
\newcommand{\econ}{\end{conjecture}}
\newcommand{\bcor}{\begin{corollary}}
\newcommand{\ecor}{\end{corollary}}
\newcommand{\bdf}{\begin{definition}}
\newcommand{\edf}{\end{definition}}
\newcommand{\beq}{\begin{equation}}
\newcommand{\eeq}{\end{equation}}
\newcommand{\bexa}{\begin{example}}
\newcommand{\eexa}{\end{example}}
\newcommand{\bexe}{\begin{exercise}}
\newcommand{\eexe}{\end{exercise}}
\newcommand{\bfac}{\begin{fact}}
\newcommand{\efac}{\end{fact}}
\newcommand{\bite}{\begin{itemize}}
\newcommand{\eite}{\end{itemize}}
\newcommand{\blem}{\begin{lemma}}
\newcommand{\elem}{\end{lemma}}
\newcommand{\bprb}{\begin{problem}}
\newcommand{\eprb}{\end{problem}}
\newcommand{\bpro}{\begin{proposition}}
\newcommand{\epro}{\end{proposition}}

\newcommand{\bque}{\begin{question}}
\newcommand{\eque}{\end{question}}
\newcommand{\brem}{\begin{remark}}
\newcommand{\erem}{\end{remark}}
\newcommand{\bthm}{\begin{theorem}}
\newcommand{\ethm}{\end{theorem}}
\newcommand{\bmat}{\begin{matrix}}
\newcommand{\emat}{\end{matrix}}

\newcommand{\bpr}{\begin{proof}}
\newcommand{\epr}{\end{proof}}

\newcommand{\lb}{\label}

\newcommand{\comment}[1]{\,}

\newcommand{\p}{\partial}

\newcommand{\Z}{\mathbb Z}

\newcommand{\R}{\mathbb R}

\newcommand{\inft}{
\begin{tikzpicture}[scale=0.03528, baseline={([yshift=-.5ex]current bounding box.center)}]
\draw (0,-13) arc [radius=9, start angle=-45, end angle=45];
\draw (10,0) arc [radius=9, start angle=135, end angle=225];
\end{tikzpicture}}
\newcommand{\zerot}{
\begin{tikzpicture}[scale=0.03528,baseline={([yshift=-.5ex]current bounding box.center)}]
\draw (0,-4) arc [radius=9, start angle=135, end angle=45];
\draw (12,6) arc [radius=9, start angle=-45, end angle=-135];
\end{tikzpicture}}

\newcommand{\pmo}{{\pm 1}}

\title{Verification Of The Jones Unknot Conjecture Up To 24 Crossings}

\author{Robert E. Tuzun, Adam S. Sikora}

\begin{document}

\thispagestyle{empty}

\begin{abstract}
Extending upon our previous work, we verify the Jones Unknot Conjecture for all knots up to $24$ crossings. 
We describe the method of our approach and analyze the growth of the computational complexity of its different components. 
\end{abstract}

\address{244 Math Bldg, University at Buffalo, SUNY, Buffalo, NY 14260}
\email{retuzun@buffalo.edu, asikora@buffalo.edu}

\pagestyle{myheadings}

\maketitle


\section{Introduction}

The Jones Unknot conjecture states that the Jones polynomial distinguishes all non-trivial knots from the unknot. 
It is a striking statement implying that the Jones polynomial -- which is a quantum-physics inspired and still mysterious invariant of knots -- contains very strong topological information. Jones proposed it as one of the challenges for mathematics in the 21st century in \cite{Jo}.

The main result of this paper is:

\begin{theorem}
The Jones polynomial distinguishes all non-trivial knots with diagrams up to 24 crossings from the unknot. 
\end{theorem}

This is an extension of our earlier work \cite{TS}, which contains broad background information and setup than the current paper. Therefore, we encourage unfamiliar reader to use it as a reference for further details on this projct. That paper includes also references to other's work on this conjecture.

\section{Description of the method}

A {\em $2$-tangle} is a $1$-manifold $T$ embedded into a $3$-ball $B^2\times [-1,1]$ with four ends at points NW, NE, SE, SW of $\p B^2\times\{ 0\}$. Tangles are considered up to isotopy in $B^3$ fixing their ends.

{\em Integral tangles} are those obtained from the zero tangle $\zerot$ by a sequence of additions of $\pm 1$, cf. Figure \ref{f-basictangles}.
Tangles obtained from integral ones by additions and $90^o$ rotations are called {\em algebraic,}  \cite{Co}.

\newcommand{\KP}[1]{%
  \begin{tikzpicture}[baseline=-\dimexpr\fontdimen22\textfont2\relax]
  #1
  \end{tikzpicture}}
  
\newcommand{\Taddition}[2]{%
  \KP{%
    \filldraw[color=black, fill=none, thick] circle (0.28);
    \filldraw[color=black, fill=none, thick] (0.8,0) circle (0.28);
    \put(-0.07,-0.05){#1}
    \put(0.235,-0.05){#2}
    \draw[color=black, thick] (-0.3,-0.3) -- (-0.19,-0.19);
    \draw[color=black, thick] (-0.3,0.3) -- (-0.19,0.19);
    \draw[color=black, thick] (1.1,-0.3) -- (.99,-0.19);
    \draw[color=black, thick] (1.1,0.3) -- (0.99,0.19);
    \draw[color=black, thick] (0.2,0.19) .. controls (0.4, 0.3)  .. (0.60,0.19);
    \draw[color=black, thick] (0.2,-0.19) .. controls (0.4,-0.3)  .. (0.60,-0.19);
  }%
}

\begin{figure}[h]
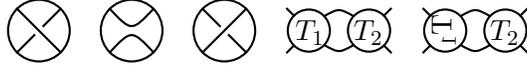

\begin{center}
\KP{
   \filldraw[color=black, fill=none, thick] (0,0) circle (0.41);
    \draw[color=black, thick] (-0.3,0.3) -- (0.3,-0.3);
    \draw[color=black, thick] (-0.3,-0.3) -- (-0.05,-0.05);
    \draw[color=black, thick] (0.05,0.05) -- (0.3,0.3);
  }\quad
  \KP{
   \filldraw[color=black, fill=none, thick] (0,0) circle (0.41);
    \draw[color=black, thick] (-0.3,0.3) .. controls (0,-0.02) .. (0.3,0.3);
    \draw[color=black, thick] (-0.3,-0.3) .. controls (0,0.02) .. (0.3,-0.3);
  }\quad
  \KP{
  \filldraw[color=black, fill=none, thick] (0,0) circle (0.41);
  \draw[color=black, thick] (-0.3,-0.3) -- (0.3,0.3);
   \draw[color=black, thick] (-0.3,0.3) -- (-0.05,0.05);
  \draw[color=black, thick] (0.05,-0.05) -- (0.3,-0.3);
  }\quad
  \Taddition{$T_1$}{$T_2$}\quad
  \Taddition{\raisebox{-.1cm}{\rotatebox{90}{\reflectbox{$\rm{T}_1$}}}}{$T_2$}
  \end{center}
\caption{The $-1, 0, 1$, 
tangles, the tangle addition, $T_1+T_2$, and multiplication, $T_1\cdot T_2$}
\lb{f-basictangles}
\end{figure}
Algebraic tangles are closed under the tangle multiplication, $T_1\cdot T_2$, which is the result of an addition of 
\bite
\item  a reflection of $T_1$ with respect of the NW-SE diagonal, and of
\item $T_2,$ 
\eite
cf. Figure \ref{f-basictangles}(right).

A closure of an algebraic tangle is an {\em algebraic} or {\em arborescent link}.

A {\em Conway polyhedron} is an edge-connected 4-valent simple planar graph with no regions with just two vertices. 
(Such polyhedra are called simple polyhedra in \cite{Co}.)
An important observation of Conway is that each knot is either algebraic or has a diagram obtained by filling the 4-valent vertices of a Conway polyhedron with algebraic tangles. Roughly speaking we use that approach to generate all necessary knot diagrams. A Conway polyhedron is {\em thin} if if it can
be disconnected by removing two of its edges. Since, without loss of generality one can consider the Jones Unknot Conjecture for prime knot diagrams only, we restrict our attention to non-thin Conway polyhedra only.

Recall that the Kauffman bracket $[\cdot ]$ is a $\Z[A^{\pm 1}]$-valued invariant of framed unoriented links in $\R^3$ defined uniquely by the following skein relations:
$$\left[\, \KP{
    \draw[color=black, thick] (-0.3,0.3) -- (0.3,-0.3);
    \draw[color=black, thick] (-0.3,-0.3) -- (-0.05,-0.05);
    \draw[color=black, thick] (0.05,0.05) -- (0.3,0.3);
  } \, \right]=
  A\left[\, \KP{
    \draw[color=black, thick] (-0.3,0.3) .. controls (0,-0.02) .. (0.3,0.3);
    \draw[color=black, thick] (-0.3,-0.3) .. controls (0,0.02) .. (0.3,-0.3);
  }\, \right]\quad
  +A^{-1}\left[\,  \KP{
  \draw[color=black, thick] (-0.3,-0.3) -- (0.3,0.3);
   \draw[color=black, thick] (-0.3,0.3) -- (-0.05,0.05);
  \draw[color=black, thick] (0.05,-0.05) -- (0.3,-0.3);
  }\, \right],\quad [L\cup \bigcirc] = (-A^2+A^{-2})[ L]
  $$ 

The notion of the Kauffman bracket can be extended to framed tangles, through the concept of the Kauffman bracket skein module of a thickened disk with four boundary marked points. It consists of the space of formal
$\Z[A^{\pm 1}]$-linear combinations of tangle diagrams in $B^2$ (with endpoints at NW, NE, SE, SW), subject to regular isotopy and the above Kauffman bracket skein relations. For a further discussion of this skein module see  for example \cite{SW}, where it is denoted by $\mathbb A_1(D^2,\{NW, NE, SE, SW\},\Z[A^{\pm 1}])$. 

By \cite[Cor. 4.1]{SW}, this $\Z[A^{\pm 1}])$-module is free with a basis given by the zero and infinity tangles, \zerot\ and \inft. In other words,  every tangle, considered as an element of the skein module,
can be written as
$$T=p\cdot \zerot+ q\cdot \inft$$
where $p,q\in \Z[A^{\pmo}]$ are uniquely defined.
We call $[T]=(p,q)$ the {\em Kauffman bracket} of $T.$

A tangle $T$ is {\em algebraically trivializable} if its Kauffman bracket $[T]=(p(A),q(A))$ is such that 
$p(A)r(A) + q(A)s(A) = 1$ for some $r(A),s(A)\in \Z[A^{\pm 1}]$.
For example, $T$ with $[T]=(3A, 2)$ is algebraically trivializable since $3A \cdot A^{-1} - 2 \cdot 1 = 1$,
but $T$ with $[T]=(4A, 2)$ is not. 

The importance of this  property stems from the fact that any potential counterexample to Jones Unknot conjecture
is either a closure of an algebraically trivializable algebraic tangle or can be obtained by filling the vertices of a Conway polyhedron with algebraically trivializable algebraic tangles. Additionally, we can assume that these tangles have no internal loops, like for example the (2,2)-pretzel tangle, since such tangles cannot be completed to a knot.

The method used in this work was that of \cite{TS} with further optimizations implemented. (See 
\cite{TS} for further details.)
 It involved the following steps:
\begin{enumerate} 
\item Generation of non-thin Conway polyhedra up to 24 vertices.
\item Generation of algebraically-trivializable algebraic tangles, up to 18 crossings, cf. Table 2. (Larger tangles, up to 24 crossings, are generated on the fly.)
\item Generation of knot diagrams of 24 crossings, by
\begin{enumerate}  
\item considering all possible insertions of algebraically-trivializable algebraic tangles into Conway polyhedra resulting in 24 crossing diagrams and by 
\item considering closures of all 24-crossing algebraic tangles. 
\end{enumerate}
\item Elimination of knot diagrams allowing a pass move which reduces either the number of crossings or the number of vertices in the corresponding Conway polyhedron.
\item Computation of the determinants of the remaining knot diagrams by a divide-and-conquer method.
\item Computation of the Kauffman bracket polynomials of determinant $1$ knot diagrams by a similar divide-and-conquer method. The diagrams with monomial Kauffman brackets are called {\em candidates}.
\item  Computation of the knot group presentations for the candidates $K$ using the computer program SnapPy \cite{Sn}. This program is higly efficient in finding the cyclic $\langle a \rangle$ presentation of $\pi_1(S^3-K)$, thus confirming the triviality of $K.$ (This step is repeated if necessary, as the SnapPy algorithm is non-deterministic.)
\item Using other methods for knot diagrams for which SnapPy did not find the cyclic presentation of their knot groups, including Dynnikov's unknot recognition algorithm, \cite{Dy}.
\end{enumerate}

Among the optimizations introduced for the 23 and 24 crossing knot testing was the analysis of mutations among Conway polyhedra and elimination of those related by mutation to other polyhedra already on the list and related to $4$-valent graphs which contain a bigon.

\section{Computational Complexity}

Computational complexity depends in part on the growth of the number of Conway polyhedra with the number of vertices (roughly between $3^{v}$ and $4^{v}$, as indicated in Table 1),
and the number of trivializable algebraic tangles with the number of crossings (roughly about
$3^{n}$, as indicated in Table \ref{algtang_table}).  It should be noted that algebraic tangles with internal
loops, such as the (2,2)-pretzel tangle, are not allowed since they would result in
multi-component links.

\begin{table}[H]
\label{polyh_table}
\title{Table 1.  Numbers of non-thin Conway polyhedra with $v$ vertices\\.}

{\begin{tabular}{|rr|rr|}
\hline
$v$ & Total & $v$ & Total \\ \hline
 6 &      1 & 16 &     499 \\
 8 &      1 & 17 &    1473 \\
 9 &      1 & 18 &    4974 \\
10 &      3 & 19 &   16296 \\
11 &      3 & 20 &   56102 \\
12 &     12 & 21 &  192899 \\
13 &     19 & 22 &  674678 \\
14 &     63 & 23 & 2381395 \\
15 &    153 & 24 & 8468424 \\
\hline
\end{tabular}}
\end{table}

\begin{table}[H]
\label{algtang_table}
\title{Table 2.  Numbers of algebraic tangles (total and trivializable) with $c$ crossings (and no internal loops).}

{\begin{tabular}{|rrr|rrr|}
\hline
$c$ & Total & Triv & $c$ & Total & Triv \\
\hline
1  &         1 &       1 & 10 &      4334 &     2589 \\
2  &         1 &       1 & 11 &     15076 &     7754 \\
3  &         2 &       2 & 12 &     53648 &    23572 \\
4  &         4 &       4 & 13 &    193029 &    71124 \\
5  &        12 &      12 & 14 &    698590 &   211562 \\
6  &        36 &      30 & 15 &   2560119 &   633059 \\
7  &       113 &      94 & 16 &   9422500 &  1866458 \\
8  &       374 &     288 & 17 &  34935283 &  5478404 \\
9  &      1242 &     836 & 18 & 130250565 & 15674910 \\ \hline
\end{tabular}}
\end{table}

\begin{figure}[H]
\includegraphics[width=0.90\textwidth]{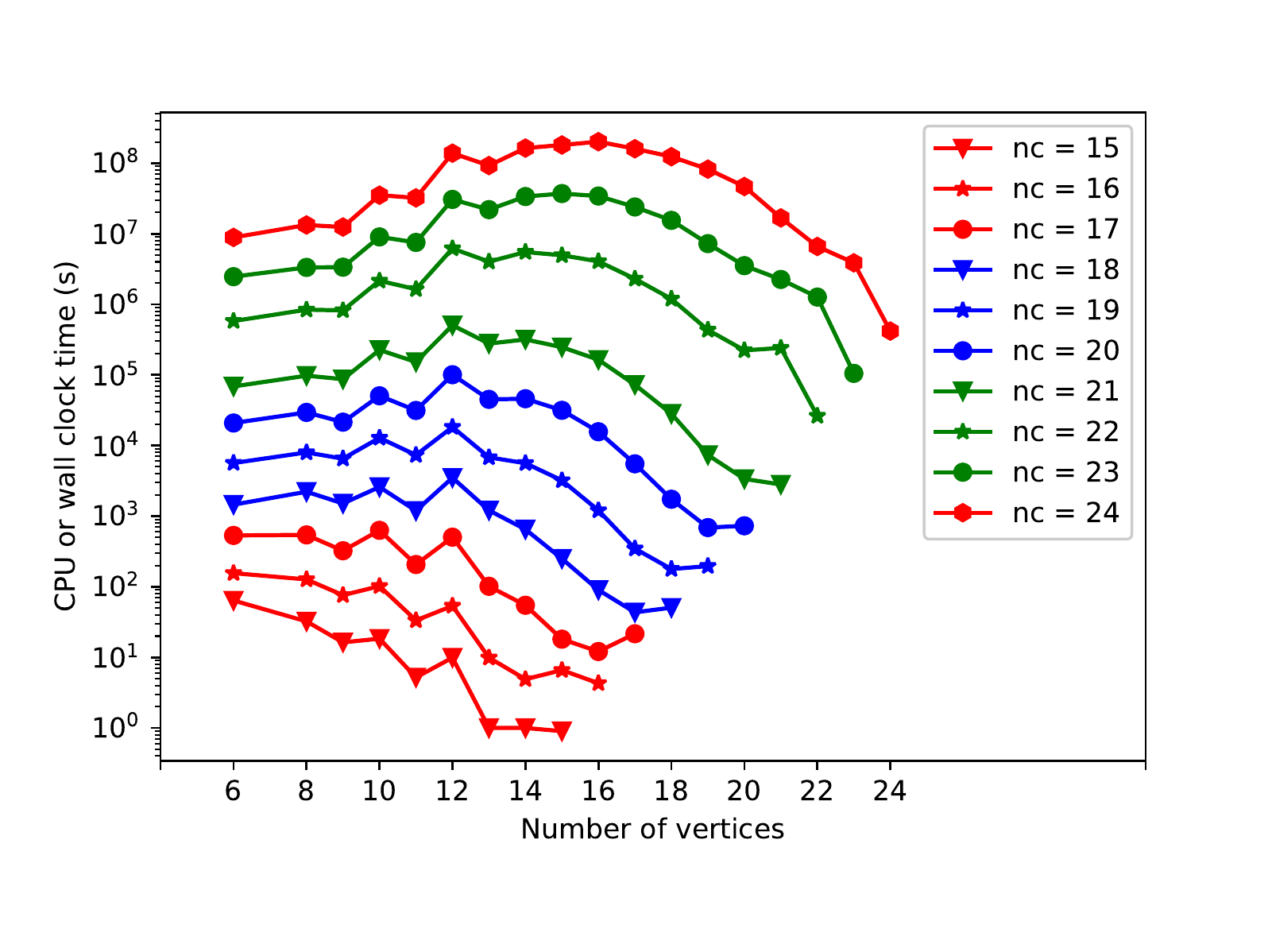} \vspace*{-.1in}\\
\caption{CPU times for non-algebraic knot diagrams generation and determinant testing. ``nc'' is the total crossing number.}
\label{times}
\end{figure}

\begin{figure}[h]
\includegraphics[width=0.90\textwidth]{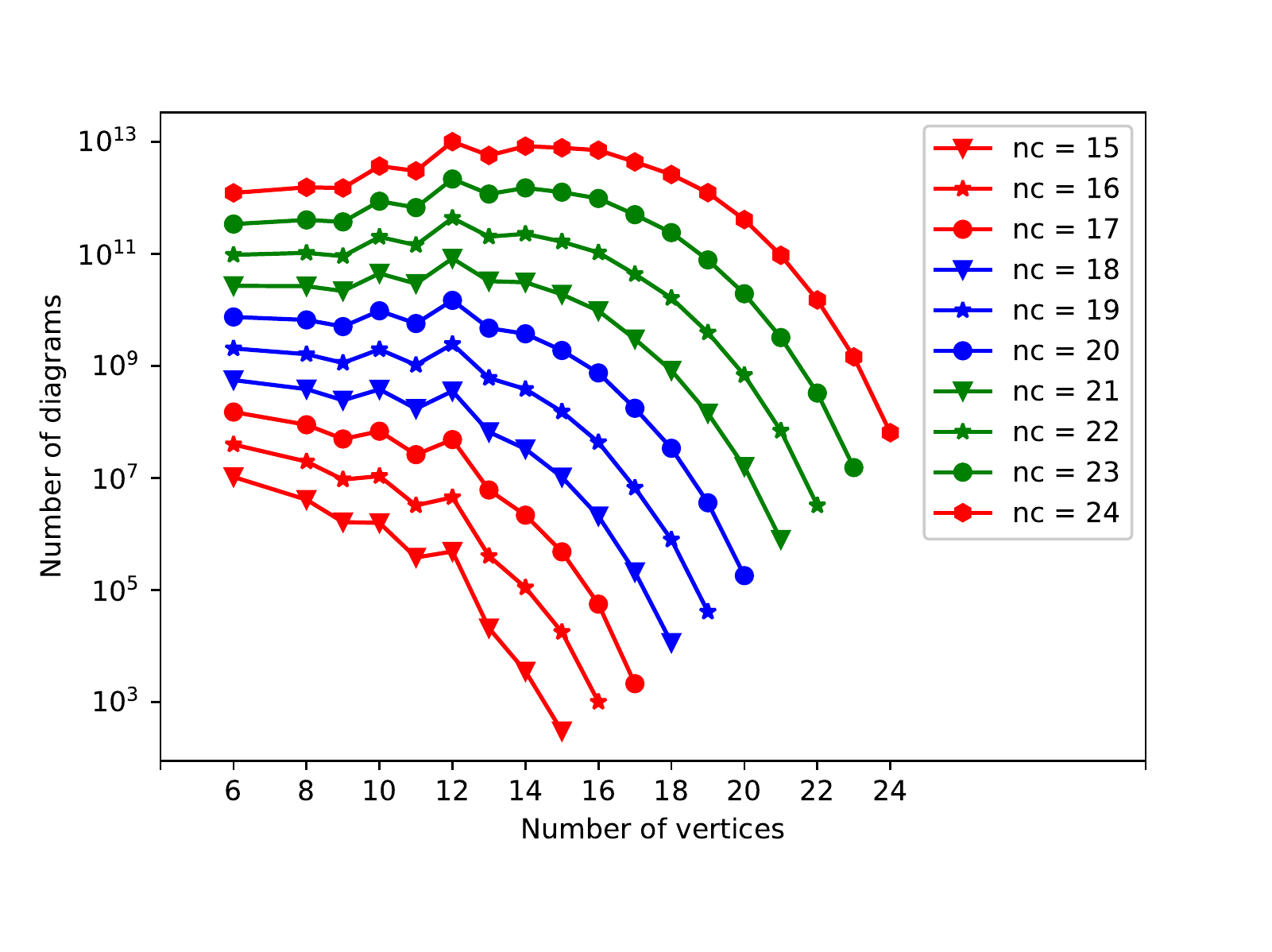}
\caption{Number of non-algebraic knot diagrams tested for determinant.}
\label{num_tested}
\end{figure}

\begin{figure}[H]
\includegraphics[width=0.95\textwidth]{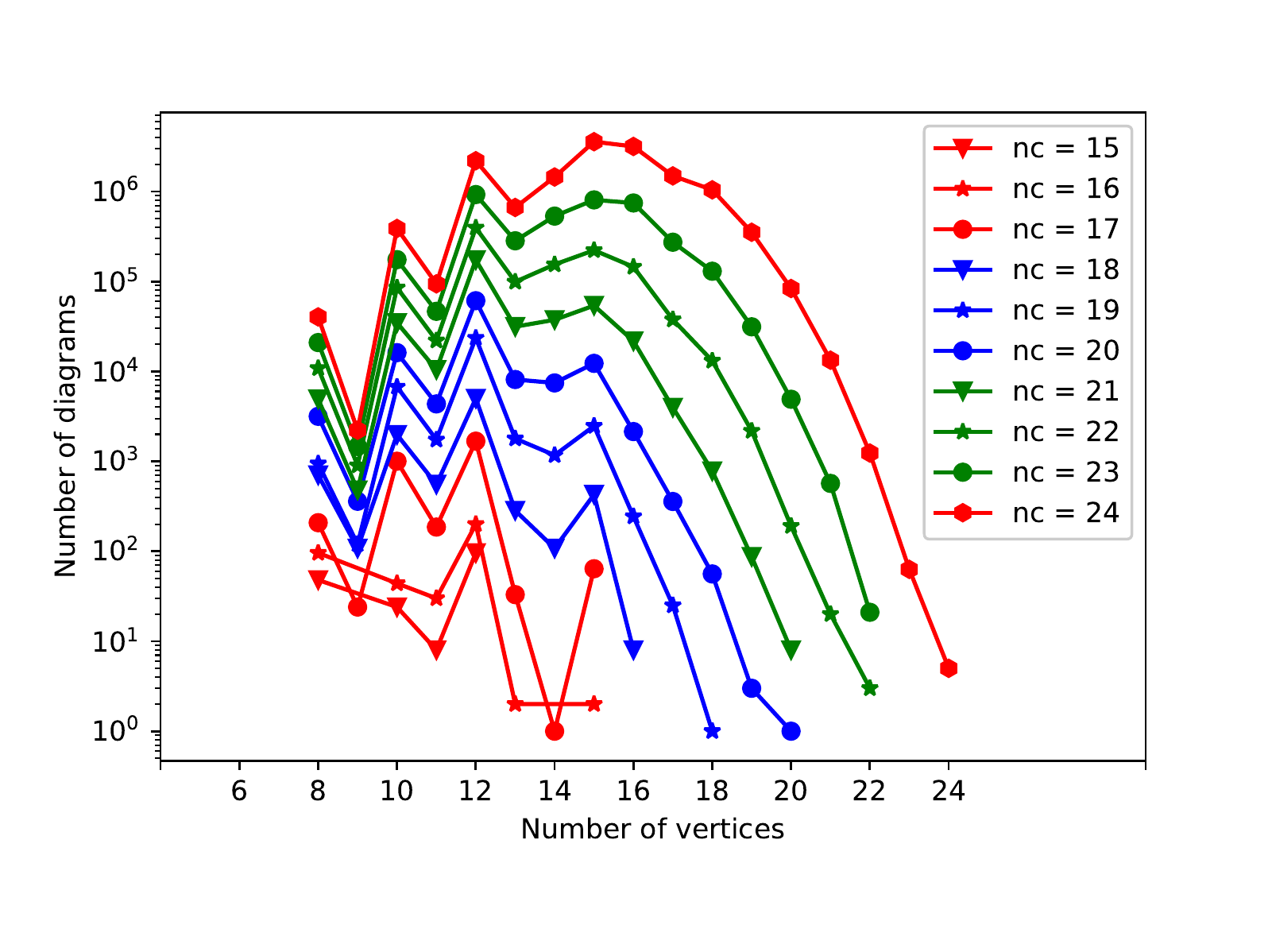}
\caption{Number of diagrams with monomial Kauffman bracket tested for unknottedness.}
\label{num_cand}
\end{figure}

Although the number of $v$ vertex Conway polyhedra grows with $v$, Conway polyhedra with the larger number of vertices admit fewer knot diagrams with a given crossing number.
The overall effect of these two patterns on determinant computation times is shown in Figure
\ref{times}.  The sudden jump in the computation times between 21 and 22 crossings
appears to be due to the relative speed of two different CPU types (i7-4790 for
calculations for 21 or fewer crossings, and Intel Xeon L5520 for 22 or more crossings).
The overall computation time for a $c$ crossing calculation is exponential in $c$, somewhere between
$5^{c}$ and $7^{c}$.  Similar trends occur for number of diagrams tested (Figure
\ref{num_tested}) and number of diagrams with monomial Kauffman brackets (Figure
\ref{num_cand}). 

\section{Summary}

The verification for 24 crossings required testing 59,361,435,729,041 non-algebraic knot diagrams and 185,317,928,640 algebraic knot diagrams, for a total of
59,546,753,657,681 knot diagrams -- approximately 6 times the number of knot
diagrams for 23 crossings.

Computations for 24 crossings were performed on 8-core Intel Xeon L5520 processors operated by the Center for Computational Research at the University at Buffalo for a total of 41.8 core-years of wall-clock time.

%

\end{document}